\documentclass{article}

\usepackage{amssymb,amsmath,makeidx,fancybox,color,colordvi,multicol,graphicx,wrapfig,soul}
\usepackage{cite}
\usepackage{standalone}
\usepackage{float}
\usepackage{subcaption}
\usepackage{amsfonts}

\let\temp\phi
\let\phi\varphi
\let\varphi\temp

\let\temp\epsilon
\let\epsilon\varepsilon
\let\varepsilon\epsilon

\newcommand{\R}{\mathbb{R}}
\newcommand{\C}{\mathbb{C}}

\newtheorem{theorem}{Theorem}[section]

\newtheorem{lemma}[theorem]{Lemma}

\newtheorem{remark-definition}[theorem]{Remark and Definition}
\title{Equivariant Versions of Odd Number Theorem}

\author{Edward Hooton\thanks{Department of Mathematical Sciences, The University of Texas at Dallas, TX, USA}
	\and Pavel Kravetc$^*$ \and Dmitrii Rachinskii$^*$}

\date{}

\begin{document}

\maketitle

	\begin{abstract}
		We consider the problem of stabilization of unstable periodic solutions to autonomous systems by the non-invasive delayed feedback control known as Pyragas control method.
		The Odd Number Theorem imposes an important restriction upon the choice of the gain matrix by stating a necessary condition for stabilization. 
		In this paper, the Odd Number Theorem is extended to equivariant systems.
		We assume that both the uncontrolled and controlled systems respect a group of symmetries. Two types of results are discussed.
		First, we consider rotationally symmetric systems for which the control stabilizes the whole orbit of relative periodic solutions that form an invariant two-dimensional torus in the phase space.
		Second, we consider a modification of the Pyragas control method that has been recently proposed for systems with a finite symmetry group. This control acts non-invasively on one selected periodic solution from the orbit and targets to stabilize this particular solution. Variants of the Odd Number Limitation Theorem are proposed for both above types of systems. The results are illustrated
		with examples that have been previously studied in the literature on Pyragas control including a system of two symmetrically coupled Stewart-Landau oscillators and a system of two coupled lasers.
	\end{abstract}
	
{\bf Keywords:} Stabilization of periodic orbits, Pyragas control, delayed feedback, $S^1$-equivariance, finite symmetry group.

\section{Introduction}

Stabilization of unstable periodic solutions is an important problem in applied nonlinear sciences.
An elegant method suggested by Pyragas \cite{pyragas1992continuous} is to introduce delayed feedback with the delay equal, or close, to the period
$T$ of the target unstable periodic solution $x^*(t)$ to the uncontrolled system $\dot x(t)=f(t,x(t))$.
This feedback control is typically linear, and the controlled system has the form
\begin{equation}\label{clas_pyr}
\dot x(t) = f(t,x(t))+ K(x(t-T)-x(t)), \qquad x\in \mathbb R^n,
\end{equation}
where  $K$ is an $n\times n$ gain matrix. Since the explicit form of the target cycle is not required this method is easy to implement and widely applicable \cite{app1,app2, app3}.
Pyragas control is often referred to as {\em non-invasive},
since $x^*(t)$ is an exact solution of both the uncontrolled and controlled systems
if the delay exactly equals the period of $x^*$. The question is how to choose the gain matrix $K$ to ensure that
$x^*$ is a stable solution of \eqref{clas_pyr}.

Certain limitations to the method of Pyragas are known.
It was proved in \cite{original_odd} that if $f$ depends explicitly on $t$ and
the target periodic solution $x^*$ of the uncontrolled non-autonomous system
is hyperbolic with an odd number of real Floquet multipliers greater than one, then for any choice of $K$,
the function $x^*$ is an {\em unstable} solution of \eqref{clas_pyr}.
In \cite{hooton_amann}, this theorem was modified to deal with the case of autonomous systems
\begin{equation}\label{1}
\dot x(t) = f(x(t))+ K(x(t-T)-x(t)).
\end{equation}
In this case, the theorem provides necessary conditions on the control matrix $K$ to allow stabilization
of an unstable hyperbolic cycle $x^*$ of the autonomous system
\begin{equation}\label{uncont1}
\dot{x}(t)=f(x(t)).
\end{equation}
These necessary conditions can be used as a guide when constructing the gain matrix $K$.

In this paper, we consider systems \eqref{uncont1}, which respect some symmetry.
Periodic solutions (cycles) of such systems naturally come in group theoretic orbits, hence there are multiple cycles with the same period.
This complicates the applicability of Pyragas control because the control acts non-invasively on all those cycles.
In particular, for systems with a continuous group of symmetries, a connected continuum of cycles, which all have the same period, is generic.
The cycles that form the continuum are not hyperbolic and hence do not satisfy
conditions of theorems from \cite{original_odd, hooton_amann}.

On the other hand, a modification of Pyragas control was proposed in \cite{fiedler_Z_2} for systems with a finite group of symmetries
in order to make the control non-invasive only on one selected target cycle, which has been chosen for stabilization.
The symmetry of a cycle $x^*$ is described by a collection of pairs $(A_g,T_g)$ where $A_g \in \text{GL}(n)$ and $T_g$ is a rational fraction of the period of $x^*$. The symmetry is expressed by the property that
\begin{equation}\label{sym_prop2}
A_gx^*(t)= x^*(t+T_g)
\end{equation}
for all the pairs $(A_g,T_g)$.
To stabilize $x^*$, it was suggested in \cite{Isabelle_thesis} to modify \eqref{1} by selecting one particular $g$
and introducing control as follows:
\begin{equation}\label{fied_cont}
\dot x(t) = f(x(t))+ K(A_gx(t-T_g)-x(t)).
\end{equation}
In \cite{fiedler_Z_2,mamanya1, mamanya_physica_d} this type of control was applied to stabilize small amplitude cycles born via a Hopf bifurcation, while in \cite{baryshnya_book_chapter} analysis of the stability of large amplitude cycles was done by exploiting the additional $S^1$ symmetry of Stewart-Landau oscillators.

In this paper, we extend the odd-number limitation type results considered in \cite{hooton_amann} to treat the case when control of the form \eqref{fied_cont}
is applied to a system with a finite group of symmetries (Section 2);
and, the case when the standard Pyragas control such as in \eqref{1}
is applied to a target cycle, which is not hyperbolic, because the system is $S^1$-symmetric (Section 3).
Analytical results are illustrated with examples.

\section{Modified Pyragas control of systems with finite symmetry group}
\subsection{Necessary condition for stabilization}
Suppose that system \eqref{uncont} has a periodic solution $x^*$ with period $T$. Assume that this system respects some group of symmetries, and for one particular $g$
relation \eqref{sym_prop2} holds.
We denote by $\Phi(t)$ the fundamental matrix of the linearization
\begin{equation}\label{lin}
\dot{y}=B(t)y,\qquad B(t):= D f_x(x^*(t)),
\end{equation}
of system \eqref{uncont} near $x^*(t)$, where $D f_x$ denotes the Jacobi matrix of $f$.
Condition \eqref{sym_prop2} implies that
\begin{equation}\label{cor}
A_g^{-1}\Phi(T_g) \psi(0)=\psi(0),\qquad \psi(t):=\dot x^*(t),
\end{equation}
i.e.~the matrix $A_g^{-1}\Phi(T_g)$ has an eigenvalue $1$. We assume that

\bigskip
$(H_1)$ \ {\em $1$ is a simple eigenvalue for the matrix $A_g^{-1}\Phi(T_g)$.}

\bigskip
Following \cite{schneider_Z_3}, we introduce a modified Pyragas control as in \eqref{fied_cont},
where we assume that
\begin{equation}\label{commute}
KA_g=A_gK.
\end{equation}
This commutativity property can be a natural restriction on feasible controls.
For example, it is typical of laser systems.
On the other hand, gain matrices, which are simple enough to allow for efficient analysis of stability
of the controlled equation \eqref{fied_cont}, also usually satisfy condition \eqref{commute} (cf.~\cite{mamanya1, schneider_elim}).

Let $D'$ denote the transpose of a matrix $D$.
Using $(H_1)$, denote by $\psi^\dagger$ the normalized adjoint eigenvector with the eigenvalue $1$ for the matrix $[A_g^{-1}\Phi(T_g)]'$:
\[
[A_g^{-1}\Phi(T_g)]' \psi_0^\dagger=\psi_0^\dagger;\qquad  \psi_0^\dagger \cdot \psi(0) =1,
\]
where dot denotes the standard scalar product in $\mathbb{R}^n$.
Furthermore, denote by $\psi^\dagger(t)$ the solution of the initial value problem
\[
\dot y= -B'(t) y; \qquad y(0)=\psi_0^\dagger.
\]
Since the fundamental matrix of system $\dot y= -B'(t) y$ is $(\Phi^{-1}(t))'$,
\begin{equation}\label{phi'}
\psi^\dagger(t)=(\Phi^{-1}(t))'\psi^\dagger_0.
\end{equation}
Note that relation \eqref{sym_prop2} implies
\[
A_g \psi(t)= \psi(t+T_g);\qquad A_g'\psi^\dagger(t)= \psi^\dagger(t-T_g).
\]

Finally denote by $N$ the number of real eigenvalues $\mu$ of the matrix $A_g^{-1}\Phi(T_g)$, which satisfy $\mu>1$.

\begin{theorem}\label{t1}
Assume that conditions $(H_1)$ and \eqref{commute} hold. Let
\begin{equation}\label{main}
(-1)^N\left(1 + \int_{0}^{T_g}  \psi^\dagger (t)\cdot K \psi(t)\,dt\right)<0.
\end{equation}
Then, $x^*(t)$ is an {\em unstable} periodic solution of the controlled system \eqref{fied_cont}.
\end{theorem}
Hence, the inequality opposite to \eqref{main} is a necessary condition for stabilization of the periodic solution $x^*$.
This necessary condition restricts the choice of the gain matrix $K$.

\subsection{Example}\label{example}
As an illustrative example  of Theorem \ref{t1} we consider the system of two identical diffusely coupled Landau oscillators described in complex form by 
\begin{equation}
\label{uncont}
\begin{aligned}
\dot z_1 &= (\alpha +i + \gamma|z_1|^2)z_1 + a(z_2 -z_1),\\
\dot z_2 &= (\alpha +i + \gamma|z_2|^2)z_2 + a(z_1 -z_2)
\end{aligned}
\end{equation}
with $z_1,z_2\in \mathbb{C}$. Here $\alpha$ and $a>0$ are real parameters while $\gamma$ is a complex parameter with ${\rm Re}\,\gamma>0$. When $\alpha$ is treated as a bifurcation parameter, this system undergoes two sub-critical Hopf bifurcations, the first at $\alpha = 0$ giving rise to a fully synchronized branch of solutions and the second at $\alpha = 2a$ giving rise to an anti-phase branch. The anti-phase branch is defined for $\alpha<2a$ and is given explicitly by the formula
\begin{equation}\label{tt}
z_1^*(t) = - z_2^*(t)=r(\alpha)e^{i\omega(\alpha) t},
\end{equation}
where $r=r (\alpha):= \sqrt{({2a-\alpha})/{{\rm Re}\, \gamma}}$ and $\omega=\omega (\alpha):= 1+r^2(\alpha) {\rm Im}\,\gamma$.  
In \cite{fiedler_Z_2} this branch was stabilized by introducing equivariant Pyragas control to system \eqref{uncont} in the following way:
\begin{equation}
\label{cont}
\begin{aligned}
\dot z_1 &= (\alpha +i + \gamma|z_1|^2)z_1 + a(z_2 -z_1) + b(z_2(t-{\pi}/{\omega})-z_1)),\\
\dot z_2 &= (\alpha +i + \gamma|z_2|^2)z_2 + a(z_1 -z_2) + b(z_1(t-{\pi}/{\omega})-z_2))
\end{aligned}
\end{equation}
with a complex gain parameter $b$.
To study stability of the anti-phase branch close to the bifurcation point in system \eqref{cont}, linear stability analysis of the origin combined with explicit knowledge of the branch made it possible to find sufficient conditions under which for some interval of $\alpha$ sufficiently close to $2a$ the branch is stable.

Due to the simple nature of the Landau oscillator, explicit computation of the fundamental matrix of the linearization of system \eqref{uncont} near the solution \eqref{tt} allows us to compute expression \eqref{main}. This gives us the following  necessary condition for stabilization of any of the anti-phase cycles \eqref{tt}: 
\begin{equation}
\label{cond}
1+\pi\,\frac{{\rm Re} \, \gamma \, {\rm Re}\, b+ {\rm Im}\, \gamma \, {\rm Im}\, b}{ \omega{\rm Re }\, \gamma}  < 0.
\end{equation}

At the point $\alpha = 2a$, $\omega = 1$, this coincides with  formula (6.10) from \cite{fiedler_Z_2} which, together with other conditions, defines the stability domain for small periodic orbits. In particular, for a fixed gain parameter $b$,
condition \eqref{cond} provides an upper bound for the interval of $\alpha$ where the anti-phase cycle is stable.

\subsection{Proof of Theorem \ref{t1}}

Linearizing system \eqref{fied_cont} near $x^*$ gives
\begin{equation}\label{fied_cont_lin}
\dot y(t) = B(t)y(t)+ K(A_gy(t-T_g)-y(t)).
\end{equation}
To prove that $x^*$ is an unstable periodic solution of \eqref{fied_cont} we will show that system \eqref{fied_cont_lin} has a solution
\begin{equation}\label{floq_type}
y_\mu^*(t) = \mu^{t/T_g}p(t),\qquad A_g p(t-T_g)=p(t)
\end{equation} 
with $\mu>1$, where the relation $A_g p(t-T_g)=p(t)$ ensures that $p$ is periodic.
It is easy to see that if the ordinary differential system
\begin{equation}\label{ode_eq}
\dot y = \left(B(t)+(\mu^{-1}- 1)K \right) y
\end{equation}
has a solution $y_\mu$ of type \eqref{floq_type}, then $y_\mu$ is also a solution of \eqref{fied_cont_lin}.
Denote by $\Psi_\mu(t)$ the fundamental matrix of \eqref{ode_eq}.

\begin{lemma}\label{l1}
If for some $\mu>1$ the matrix $ A_g^{-1}\Psi_\mu (T_g)$ has the eigenvalue $\mu$, then system \eqref{ode_eq} has a solution of type \eqref{floq_type} and hence 
the periodic solution $x^*$ of \eqref{fied_cont} is unstable.
\end{lemma}

The proof of this lemma is presented in the next subsection.
In order to use Lemma \ref{l1}, we consider the characteristic
polynomial 
$$
F(\mu):=\det\left(\mu\, \mathbb{I}\text{d} - A_g^{-1}\Psi_\mu(T_g)\right)
$$
of the matrix $ A_g^{-1}\Psi_\mu (T_g)$. Observe that 
equation \eqref{ode_eq} with $\mu=1$ coincides with \eqref{lin},
hence $\Psi_1=\Phi$ and therefore condition ($H_1$) implies $F(1)=0$.
We are going to show that relation \eqref{main} implies
\begin{equation}\label{<}
F(1+\varepsilon)<0,\qquad 0<\varepsilon\ll 1.
\end{equation}
Since  $F(\mu) \to +\infty$ as ${\mu \to +\infty}$, relation \eqref{<} implies that $F$ has a root $\mu>1$
and therefore the conclusion of Theorem \ref{t1} follows from \eqref{<} by Lemma \ref{l1}.

Setting $\mu=1+\varepsilon$ and $t=T_g$ in the identity
$$
\Psi_\mu(t)= \Phi(t)\left(\mathbb{I}\text{d} + (\mu^{-1}-1)\int_{0}^{t}  \Phi^{-1}(s)K \Psi_\mu(s)\,ds\right)
$$
and
using the fact that $\Psi_{1+\epsilon}(T_g)=\Phi(T_g)+O(\varepsilon)$, we obtain the expansion
\[
\Psi_{1+\epsilon}(T_g)=\Phi(T_g)\left(\mathbb{I}\text{d} - \varepsilon Q\right)+O(\varepsilon^2), \qquad Q:=\int_{0}^{T_g}  \Phi^{-1}(t)K \Phi(t)\,dt.
\]
Therefore,
\begin{equation}\label{asy}
\begin{array}{rcl}
F(1+\epsilon) 
=
\det\left(  \text{Id}-A_g^{-1}\Phi(T_g) +\varepsilon \left(\mathbb{I}\text{d} + A_g^{-1}\Phi(T_g) Q\right)\right)+O(\varepsilon^2).
\end{array}
\end{equation}
Let us denote by $L$ the transition matrix to a basis
in which the matrix $A_g^{-1}\Phi(T_g)$ assumes the Jordan form
and agree that $\psi(0)$ is the first vector of this basis (cf.~\eqref{cor}), i.e.
\begin{equation}\label{ps}
L e_1=\psi(0),\qquad e_1:=(1,0,\ldots,0)'\in\mathbb{R}^n.
\end{equation}
In this basis, the matrix $\text{Id}-L^{-1}A_g^{-1}\Phi(T_g)L$ has the Jordan structure 
with the diagonal entries $0, 1-\mu_2,1-\mu_3,\ldots,1-\mu_n$, where $\mu_k$ are the eigenvalues of $A_g^{-1}\Phi(T_g)$ different from the simple eigenvalue $1$.
With this notation, formula \eqref{asy} implies 
\begin{equation}\label{FF}
F(1+\epsilon)=\varepsilon M_{11} \prod^{n}_{k=2}(1-\mu_k)+O(\varepsilon^2),
\end{equation}
where
\begin{equation}\label{M}
M:=\mathbb{I}\text{d} + L^{-1} A_g^{-1}\Phi(T_g) Q L; \qquad M_{11}=e_1 \cdot Me_1. 
\end{equation}

Formula \eqref{ps} implies that $(L^{-1})'e_1=\psi_0^\dagger$, hence
\[
e_1\cdot L^{-1} Q L e_1= \int_{0}^{T_g}  \psi_0^\dagger \cdot \Phi^{-1}(t)K \Phi(t) \psi(0) \,dt.
\]
Combining this with \eqref{ps} and $\psi(t)=\Phi(t)\psi(0)$, we obtain
\[
e_1\cdot L^{-1} Q L e_1= \int_{0}^{T_g}  \psi^\dagger(t) K \psi(t) \,dt.
\]
Since the first row of the matrix $L^{-1} A_g^{-1}\Phi(T_g)L$ is $(e_1)'=(1,0,\ldots,0)$, we see from \eqref{M} that
\[
M_{11}=1+e_1\cdot L^{-1} Q L e_1= 1+\int_{0}^{T_g}  \psi^\dagger(t) K \psi(t) \,dt.
\]
Hence \eqref{FF} implies 
\[
{\rm sgn}\, F(1+\epsilon)=(-1)^N \left( 1+\int_{0}^{T_g}  \psi^\dagger(t) K \psi(t) \,dt\right),
\]
where $N$ is the number of eigenvalues $\mu_k$ satisfying $\mu_k>1$.
Thus, formula \eqref{main} indeed implies \eqref{<}. 

\subsection{Proof of Lemma \ref{l1}}


Let us denote
$$
C(t):=B(t)+(\mu^{-1}-1)K.
$$
The main ingredient for the proof is the identity $$C(t+T_g) = A_gC(t)A_g^{-1},$$
which is a simple consequence of the facts that $Df(A_gx) = A_gDf(x)A_g^{-1}$, $x^*(t+T_g)=A_gx^*(t)$, and $A_g K= KA_g$.

To complete the proof denote by $\nu$ the eigenvector of $A_g^{-1}\Psi_\mu$ with the eigenvalue $\mu$ and consider the solution $y_0(t)$ of \eqref{ode_eq} with $y_0(0)=\nu$. It is clear that $y_0(t+T_g)$ satisfies
the initial value problem
\begin{equation}
\label{timestep}
\begin{cases}
\dot y_1 = C(t+T_g)y_1= A_gC(t)A_g^{-1}y_1,\\
y_1(0)= \Psi(T_g)\nu.
\end{cases}
\end{equation}
By the change of variables $y_2=A_g^{-1}y_1$, we can see that the solution of \eqref{timestep} is given by $$y_0(t+T_g)= y_1(t)= A_g\Psi(t)A_g^{-1}\Psi(T_g)\nu.$$
However by assumption $A_g^{-1}\Psi(T_g)\nu = \mu\nu$, hence
$$
y_0(t+T_g)=\mu A_g \Psi(t)\nu = \mu A_g y_0(t),
$$
which proves the lemma.

\section{Pyragas control of systems with $S^1$ spatial symmetry}
\subsection{Necessary condition for stabilization}
Suppose that equation \eqref{uncont} is $S^1$-equivariant:
\begin{equation}\label{21}
f(e^{\theta J}x)=e^{\theta J} f(x)
\end{equation}
for all $\theta \in\mathbb{R}$, $x\in\mathbb{R}^N$, where the skew-symmetric non-zero matrix $J$ satisfies $e^{2\pi J}=\mathbb{I}{\rm d}$.
We assume that \eqref{uncont} has a periodic solution $x^*(t)$ of a period $T$, which is not a relative equilibrium. Hence, equation \eqref{uncont} has
an orbit of $T$-periodic non-stationary solutions $e^{\theta J} x^*(t+\tau)$ with arbitrary $\theta,\tau$.
Therefore, the linearization \eqref{lin} has two linearly independent zero modes (periodic solutions):
\begin{equation}\label{22}
\psi_1(t)=\dot x^*(t),\qquad \psi_2(t)=J x^*(t)
\end{equation}
with the Floquet multiplier $1$.
We additionally assume that 

\bigskip
$(H_2)$ \ {\em The eigenvalue $1$ of the monodromy matrix $\Phi(T)$ of system \eqref{lin} has multiplicity exactly $2$.}

\bigskip
Then, there are two adjoint eigenfunctions (periodic solutions of equation $\dot y =-B(t)y$) that can be normalized as follows:
\begin{equation}\label{23}
\psi_1^\dagger(t)\cdot \psi_1(t)= \psi_2^\dagger(t)\cdot\psi_2(t)\equiv 1,\qquad
\psi_1^\dagger(t)\cdot \psi_2(t)= \psi_2^\dagger(t)\cdot\psi_1(t)\equiv 0.
\end{equation}

\begin{theorem}\label{t2}
Assume that condition $(H_2)$ holds. Let
\begin{equation}\label{s1maincond}
(-1)^N \bigl( (1+c_{11})(1+ c_{22})-c_{12}c_{21}\bigr)<0,
\end{equation}
where $N$ is the number of real eigenvalues $\mu$ of the monodromy matrix $\Phi(T)$, which satisfy $\mu>1$, and
\begin{equation}\label{cdef}
c_{ij}=\int_{0}^{T}  \psi^\dagger (t)\cdot K \psi_j(t)\,dt.
\end{equation}
Then, $x^*(t)$ is an {\em unstable} periodic solution of the controlled system \eqref{clas_pyr}.
\end{theorem}

\subsection{Proof of Theorem \ref{t2}}
Up to the asymptotic expansion \eqref{asy} the proof of Theorem \ref{t2} is a modification of the proof of Theorem \ref{t1}  where $A_g$ is replaced by the identity matrix $\mathbb{I}\text{d}$ and $T_g$ is replaced by $T$. In this case, the counterpart of relation \eqref{asy} is given by
\begin{equation}\label{asy2}
F(1+\epsilon)
=
\det\left(  \mathbb{I}\text{d}-\Phi(T) +\varepsilon \left(\mathbb{I}\text{d} +\Phi(T) \int_{0}^{T}  \Phi^{-1}(t)K \Phi(t)\,dt\right)+O(\varepsilon^2)\right).
\end{equation}
We again denote by $L$ the transition matrix to a basis
in which the matrix $\Phi(T)$ assumes the Jordan form
and agree that $\psi_1(0)$ and $\psi_2(0)$ are the first and second  vector of this basis (cf.~\eqref{22}), i.e.
\begin{align*}\label{ps}
L e_1=&\psi_1(0),\qquad e_1:=(1,0,\ldots,0)'\in\mathbb{R}^n,\\
L e_2=&\psi_2(0),\qquad e_2:=(0,1,\ldots,0)'\in\mathbb{R}^n.
\end{align*}
In this basis, the matrix $\text{Id}-L^{-1}\Phi(T)L$ has the Jordan structure 
with the diagonal entries $0, 0 ,1-\mu_3,\ldots,1-\mu_n$, where $\mu_k$ are the eigenvalues of $\Phi(T)$, which are different from $1$.
With this notation, formula \eqref{asy2} implies 
\begin{equation}\label{FF}
F(1+\epsilon)=\varepsilon ^2 (M_{11}M_{22}-M_{12}M_{21}) \prod^{n}_{k=3}(1-\mu_k)+O(\varepsilon^3),
\end{equation}
where
\begin{equation}\label{M}
M:=\mathbb{I}\text{d} + L^{-1} \Phi(T) \left(\int_{0}^{T}  \Phi^{-1}(t)K \Phi(t)\,dt\right)L; \qquad M_{ij}=e_i \cdot Me_j. 
\end{equation}
The same argument as in the proof of Theorem \ref{t1} shows that 
\[
{\rm sgn}\, F(1+\epsilon)=(-1)^N \bigl( (1+c_{11})(1+c_{22})-c_{12}c_{21}\bigr),
\]
where $c_{ij}$ is defined by \eqref{cdef}. Combining this with the case of Lemma \ref{l1} where $A_g = \text{Id}$ and $T_g=T$, and the fact that $F(\mu)\to+\infty$  as $\mu\to+\infty$ completes the proof.


\subsection{Example}\label{example2}

In order to illustrate Theorem \ref{t2}, we consider a model of two coupled lasers, see for example \cite{YanchukSchneider}.
In dimensionless form, the rate equations describing this system can be written as
\begin{eqnarray}
\dot{E}_1 &=& i \delta E_1 + (1+i \alpha)N_1E_1+\eta e^{-i\phi}E_2, \label{itl1} \\
\dot{N}_1 &=& \epsilon \left[J - N_1 - (1+2 N_1) |E_1|^2\right], \label{itl2} \\ 
\dot{E}_2 &=& (1+i\alpha)N_2E_2 + \eta e^{-i\phi}E_1, \label{itl3} \\
\dot{N}_2 &=& \epsilon \left[ J - N_2 - (1+2N_2)|E_2|^2\right], \label{itl4}
\end{eqnarray}
where the complex-valued variables $E_1$, $E_2$ are optical fields and the real-valued variables $N_1$, $N_2$ are carrier densities
in two laser cavities, respectively. 
This system is symmetric under the action of the group $S^1$ of transformations
$\left(E_1,N_1,E_2,N_2\right)\to \left(e^{i\theta}E_1,N_1,e^{i\theta}E_2,N_2\right)$.  Hence, the system admits relative equilibria of the form
\begin{equation}
\left(E_1, N_1, E_2, N_2\right) = \left(a_1 e^{i \omega t}, n_1, a_2 e^{i \omega t}, n_2\right) \label{relatEq}
\end{equation} with $\omega, n_1, n_2 \in \R $ and $a_1, a_2 \in \C$.
The problem of stabilization of unstable relative equilibria for this system was considered in \cite{FiedlerYanchuk}.
System \eqref{itl1}--\eqref{itl4} can also have {\it relative periodic orbits}, {\it i.e.} solutions of the form
\begin{equation}\label{tar}
	(e^{i \omega t} E_1^*(t), N_1^*(t), e^{i \omega t} E_2^*(t), N_2^*(t)),
\end{equation}
where $E_1^*(t), N_1^*(t), E_2^*(t), N_2^*(t)$ are $T$-periodic. 
In the present section we choose a relative periodic solution as a target state for stabilization.

\begin{figure}[H]
	\centering
	\includegraphics[width = .6\columnwidth]{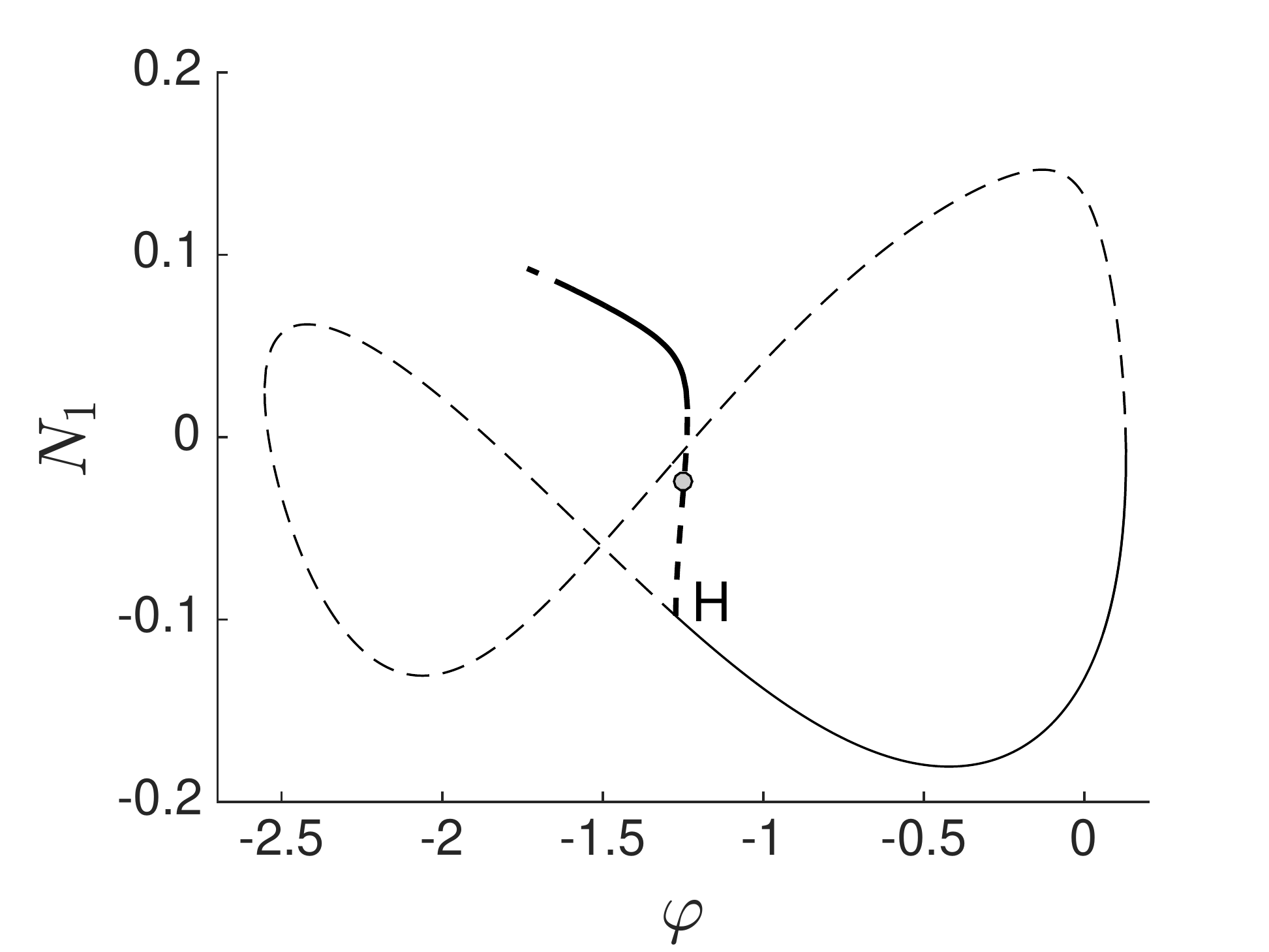}
	\caption{\small Bifurcation diagram obtained with numerical package DDE-BIFTOOL \cite{Engelborghs, 2001dde} for system \eqref{itl1}--\eqref{itl4}. Thin lines: relative equilibriua; thick lines: relative periodic solutions. Solid and dashed lines represent stable and unstable parts of the branches, respectively. H: subcritical Hopf bifurcation point; gray dot: unstable periodic orbit targeted for stabilization by Pyragas control. Parameters are $\epsilon=0.03$, $J=1$, $\eta=0.2$, $\delta=0.3$, $\alpha=2$.}
	\label{fig:itlbif}
\end{figure}

In order to stabilize the solution \eqref{tar}, we add the modified Pyragas control term
\begin{equation}  \label{quasipyragas}
	E_b(t) := b_0 \, e^{i\beta} \left( e^{-i \omega T} E_1(t-T) - E_1(t)\right)
\end{equation}
to the right hand side of equation \eqref{itl1}. Here the parameters $b_0>0$ and $\beta \in \R$ measure the amplitude and phase of the control, respectively; 
and $T$, $\omega$ are the parameters of the target relative periodic solution \eqref{tar}.
Introducing the rotating coordinates 
	$(\tilde{E}_1, \tilde{N}_1, \tilde{E}_2, \tilde{N}_2) = (e^{-i \omega t} E_1, N_1, e^{-i \omega t} E_2, N_2 )
$
transforms equations \eqref{itl1}--\eqref{itl4} 
to an autonomous system that 
has an orbit of non-stationary $T$-periodic solutions  $(e^{i\theta}\tilde{E}_1, \tilde{N}_1, e^{i\theta}\tilde{E}_2, \tilde{N}_2)=(E_1^*(t+\tau), N_1^*(t+\tau), E_2^*(t+\tau), N_2^*(t+\tau))$
with arbitrary $\theta,\tau$.
This change of variables transforms the control term \eqref{quasipyragas} to the standard Pyragas form
\begin{equation}\label{quasipyragas:rot}
\tilde{E}_b (t) = b_0 e^{i \beta} (\tilde{E}_1(t-\tau) - \tilde{E}_1(t) ).
\end{equation}
%
Hence, we can use Theorem \ref{t2} to establish the values of $b_0$ and $\beta$ for which
the control cannot stabilize the solution \eqref{tar}.

Following the analysis presented in \cite{FiedlerYanchuk}, we use the phase $\phi$ of coupling between the lasers as the bifurcation parameter. Varying $\phi$ one observes Hopf bifurcations on the branches of relative equilibria. These bifurcations give rise to branches of relative periodic solutions (which are just periodic solutions in the rotating coordinates). Figure \ref{fig:itlbif} features the bifurcation diagram for system \eqref{itl1}--\eqref{itl4} with the same parameter set as in \cite{FiedlerYanchuk}. We are interested in the unstable part of the branch of relative periodic solutions born via a subcritical Hopf bifurcation.

\begin{figure}[H]
	\centering
	\includegraphics[width=.6\columnwidth]{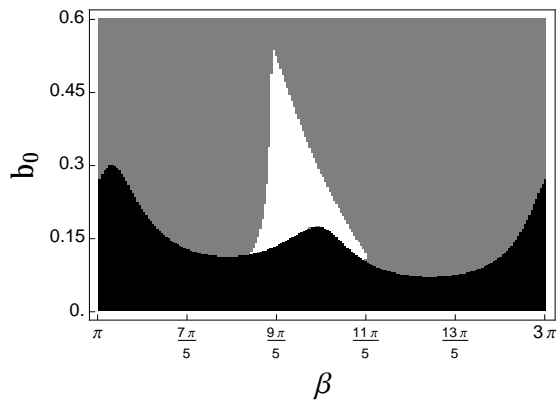}
	\caption{Domains of stability of the target relative periodic solution.
		Parameters correspond to the gray dot in Figure 1.
		Black region: sufficient condition \eqref{s1maincond} for instability is satisfied; white region: relative periodic solution is stable; gray region: relative periodic solution is unstable.}
	\label{fig:domains}
\end{figure}

\begin{figure}[H]
\begin{center}
\begin{subfigure}[b]{0.49\columnwidth}
\includegraphics*[width=\textwidth]{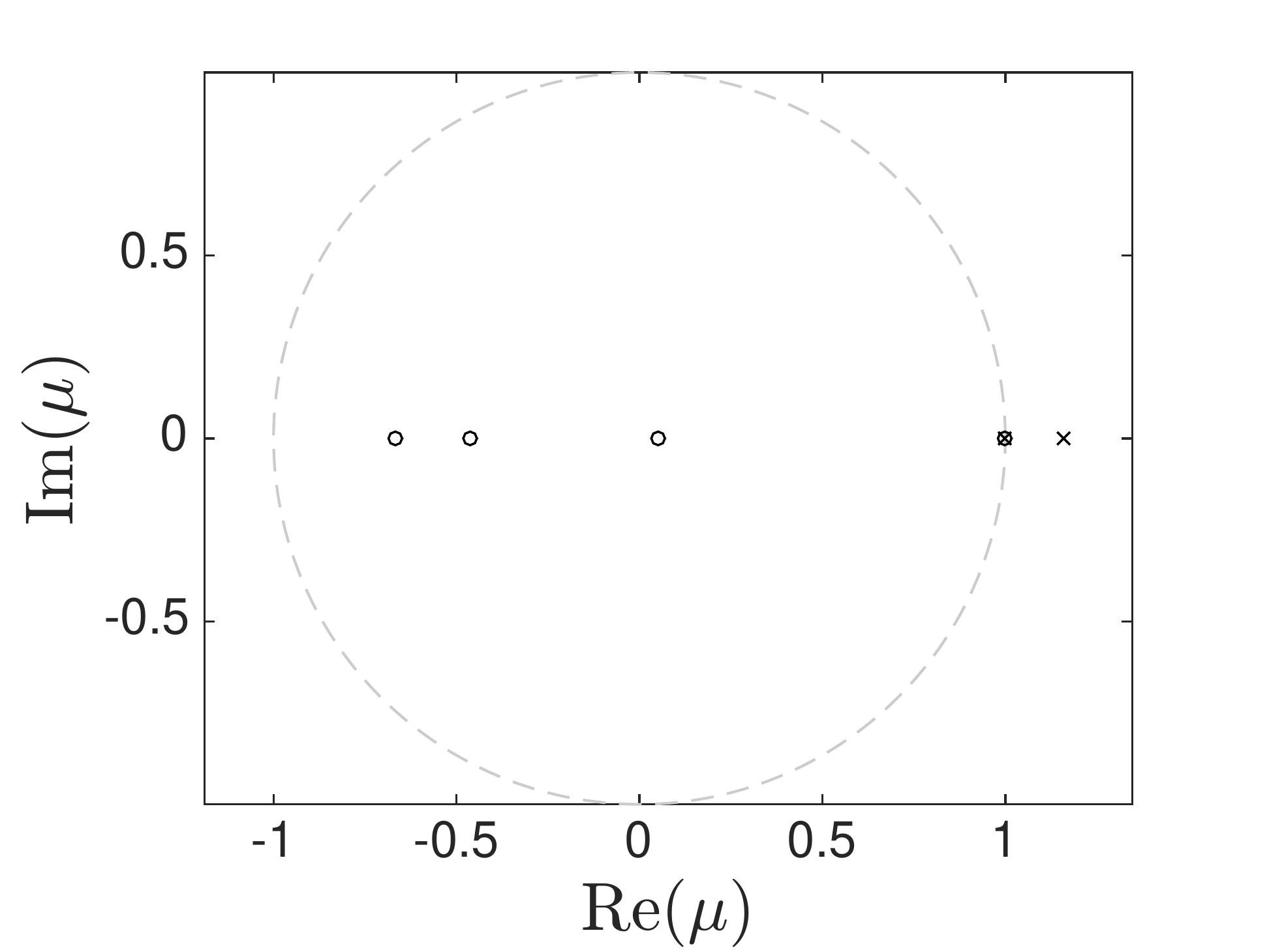}
\label{fig:unst}
\caption{}
\end{subfigure}\hfill
\begin{subfigure}[b]{0.49\columnwidth}
\includegraphics*[width=\textwidth]{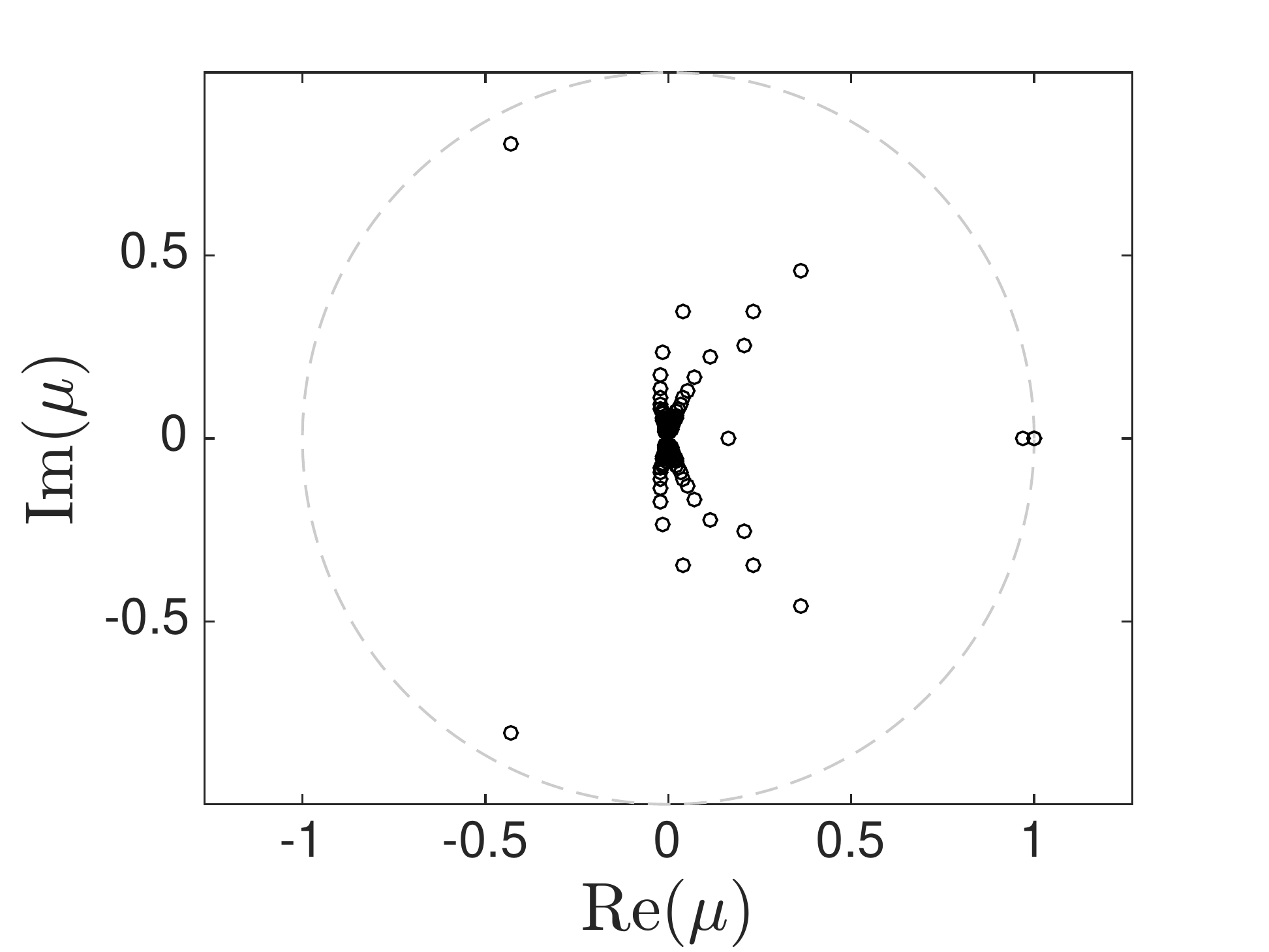}
\label{fig:st}
\caption{}
\end{subfigure}
\end{center}
\caption{Panel (a): Floquet multipliers of the target relative periodic orbit in the uncontrolled system \eqref{itl1}--\eqref{itl4}. Panel (b): Floquet multipliers of the same relative periodic orbit in the controlled system with the parameters $b_0 = 0.3036$ and $\beta = 6$ of control \eqref{quasipyragas}.}
\label{fig:s1:spectrum}
\end{figure}

Figure \ref{fig:domains} shows three regions in the parameter space $(\beta, b_0)$. The black region corresponds to the values of $b_0$ and $\beta$ for which condition \eqref{s1maincond} is satisfied, hence the target state is not stabilizable by control \eqref{quasipyragas}. The white and gray regions correspond to stable and unstable target state in the controlled system, respectively.
Figure \ref{fig:s1:spectrum} shows the change of the spectrum of the target state after applying Pyragas control \eqref{quasipyragas}.

It would be interesting to consider symmetric coupling of laser models of other types such as for example in \cite{0,1,2,3,4,5}.

\section{Discussion and Conclusions}

We have obtained necessary conditions for stabilization of unstable periodic solutions to symmetric autonomous systems by the Pyragas delayed feedback control. For systems with a finite symmetry group, we used a modification \eqref{fied_cont} of Pyragas control proposed in \cite{fiedler_Z_2}. 
This control is designed to act non-invasively on one particular solution from an orbit of periodic solutions with a specific targeted symmetry.
Inequality \eqref{main} which makes the stabilization by this control impossible is similar to its counterpart known for the standard Pyragas control
of non-symmetric systems such as \eqref{1} (see \cite{hooton_amann}):
\begin{equation}\label{main"}
(-1)^N\left(1 + \int_{0}^{T}  \psi^\dagger (t)\cdot K \psi(t)\,dt\right)<0,
\end{equation}
where $\psi=\dot x^*$ is the periodic Floquet mode of the target periodic solution $x^*$, $\psi^\dagger$ is the normalized adjoint Floquet mode, and $T$ is the period of $x^*$.
However, the necessary condition for stabilization by the modified Pyragas control (cf.~\eqref{fied_cont}) is typically more restrictive than the necessary condition for the standard Pyragas control (cf.~\eqref{1}).
This can be seen from the example of Section \ref{example}. The fully synchronized branch of cycles of system \eqref{uncont} ($z_1=z_2$), which bifurcates from zero at $\alpha=0$, can be stabilized by the standard
Pyragas control, at least near the Hopf bifurcation point \cite{fiedler_Z_2}.
The necessary condition for stabilization (i.e., the inequality opposite to \eqref{main"}) has the form
\begin{equation}
\label{cond"}
1+2\pi\,\frac{{\rm Re} \, \gamma \, {\rm Re}\, b+ {\rm Im}\, \gamma \, {\rm Im}\, b}{ \omega{\rm Re }\, \gamma}  < 0.
\end{equation}
The same control fails to stabilize the anti-phase branch of cycles ($z_1(t)=-z_2(t-T/2)$), which bifurcates from zero at another Hopf point $\alpha=2a$.
On the other hand, the modified Pyragas control (cf.~\eqref{cont}) with a proper choice of the parameter $b$ successfully stabilizes the anti-phase branch near the Hopf point \cite{fiedler_Z_2}.
Comparing the necessary conditions \eqref{cond} and \eqref{cond"} for the two controls, we see that \eqref{cond} is more restrictive because the value of the integral in \eqref{main"}
is twice the value of the integral in \eqref{main} since $T_g=T/2$ and $A_g=-\mathbb{I}{\rm d}$.
It is important to note that condition \eqref{cond} is necessary for stabilizing any cycle of the global anti-phase branch by the modified control.
In particular, it is part of the set of sufficient conditions obtained in \cite{fiedler_Z_2} for stabilizing small cycles.

We have further considered $S^1$-equivariant systems with the usual Pyragas control such as in \eqref{1}.
Due to symmetry, (relative) periodic solutions of such systems come in an $S^1$-orbit and form a two-dimensional torus in the phase space.
The control aims to stabilize all the solutions of this torus. The necessary condition for stabilization here, i.e.~the inequality opposite to \eqref{s1maincond}, is more subtle than its counterpart 
for non-symmetric systems because the solutions on the torus have the characteristic multiplier $1$ of multiplicity $2$ while in the non-symmetric case this multiplier is simple.
The control preserves the multiplier $1$ with its multiplicity.

It should be noted that any of the inequalities \eqref{main}, \eqref{s1maincond} or \eqref{main"} prevents the stabilization because it implies the existence of a real unstable characteristic multiplier $\mu>1$
(as shown in the above proofs). At the same time, our results do not help to control complex characteristic multipliers. This can be seen from the example of Section \ref{example2}.
On the border between the white stability domain and the black instability region (see Figure \ref{fig:domains}) a real characteristic multiplier passes through the value $1$, and its stability
is controlled by the sign of the left hand side of inequality \eqref{s1maincond}. On the other hand, on the border between the white domain and the gray instability region,
the change of stability is due to a pair of complex characteristic multipliers crossing the unit circle.

Stabilization of a periodic solution is more challenging in the case when the number $N$ of real characteristic multipliers which are greater than $1$ is odd than in the case when $N$ is even.
Relations \eqref{main} and \eqref{s1maincond} show that in the case of an odd $N$ the Pyragas control with a gain matrix $K$ of small norm cannot be successful, while for an even $N$ the periodic solution
may be stabilizable by small controls. On the other hand, a control with a too large amplitude is generally not successful in either case because it pushes some characteristic multipliers out of the unit circle.

\subsection*{Acknowledgments}
The authors acknowledge the support of NSF through grant DMS-$1413223$.

\bibliographystyle{plain}
\bibliography{the_bib}
\end{document}